\documentclass[12pt]{amsart}
\usepackage{amsfonts, amsthm, amsmath, amssymb, amscd, t1enc, indentfirst, graphicx, graphics, pict2e, epic, epstopdf, booktabs}
\usepackage{url}
\usepackage{enumitem}
\usepackage{algpseudocode} \algrenewcommand\algorithmicindent{0.9em}
\numberwithin{equation}{section}
\usepackage[latin2]{inputenc}
\usepackage[mathscr]{eucal}
\usepackage[margin=2.9cm]{geometry}
\usepackage{hyperref}

\theoremstyle{plain}
\newtheorem{theorem}{Theorem}[section]
\newtheorem{lemma}[theorem]{Lemma}

\theoremstyle{definition}
\newtheorem{Def}[theorem]{Definition}
\newtheorem{Conj}[theorem]{Conjecture}
\newtheorem{Rem}[theorem]{Remark}
\newtheorem{example}[theorem]{Example}
\newtheorem{notation}[theorem]{Notation}

\renewcommand{\pmod}[1]{{\ifmmode\text{\rm\ (mod~$#1$)}\else\discretionary{}{}{\hbox{ }}\rm(mod~$#1$)\fi}}

\begin{document}

\title[Algorithm for Egyptian fractions with restricted denominators]{An algorithm for Egyptian fraction representations with restricted denominators}

\author[G. Martin]{Greg Martin}
\author[Y. Shi]{Yue Shi}

\address{ Department of Mathematics \\ University of British Columbia \\ Room 121, 1984 Mathematics Road \\ Vancouver \\ BC \\ Canada \\ V6T 1Z2}
\email{gerg@math.ubc.ca}

\address{ Department of Mathematics \\ Indiana University, Bloomington \\ Bloomington \\ IN 47405 \\ United States}
\email{shi10@iu.edu}


\begin{abstract}
A unit fraction representation of a rational number~$r$ is a finite sum of reciprocals of positive integers that equals~$r$. Of particular interest is the case when all denominators in the representation are distinct, resulting in an Egyptian fraction representation of~$r$. Common algorithms for computing Egyptian fraction representations of a given rational number tend to result in extremely large denominators and cannot be adapted to restrictions on the allowed denominators. We describe an algorithm for finding all unit fraction representations of a given rational number using denominators from a given finite multiset of positive integers. The freely available algorithm, implemented in Scheme, is particularly well suited to computing dense Egyptian fraction representations, where the allowed denominators have a prescribed maximum.
\end{abstract}

\maketitle

\section{Introduction}
An {\em Egyptian fraction} is a sum of reciprocals of distinct positive integers, the name coming from the historical fact that ancient Egyptians used this form to represent rational numbers, as shown in the Rhind Papyrus~\cite{Rhindpapyrus}. Many algorithms have been developed to find Egyptian fraction representations of arbitrary positive rational numbers. Fibonacci showed in 1202 that the greedy algorithm terminates (a result rediscovered many times, most notably by Sylvester~\cite{Sylvester1880} in 1880). The equally straightforward splitting algorithm has also been proved to terminate: Beeckmans~\cite{Beeckmans1994} published a proof, although Wagon~\cite{Wagon2010} attributes the result to Graham and Jewett. An elementary pairing algorithm is also known to terminate, as proved (though not stated) by Takenouchi~\cite{Takenouchi1921}. Other algorithms include methods based on binary representations, a continued fraction algorithm, and a Farey series algorithm, along with variants; interested readers can find more information and citations in~\cite{Beeckmans1994,Bleicher1972,Eppstein1995,GN2013}.

The largest denominators in the Egyptian fraction representations resulting from these algorithms can be enormous, whereas it is often desirable to limit the size of the allowed denominators in advance. In the papers~\cite{Martin1999,Martin2000} on ``dense'' Egyptian fractions, the first author described algorithms that found Egyptian fraction representations with denominators restricted to the set $\{1,2,\dots,n\}$ for suitable positive integers~$n$. Inspired by those results, we have created a freely available implementation of such an algorithm.

Given a positive rational number~$r$ and a finite set of permitted denominators~$S$, our algorithm returns every subset of~$S$ corresponding to Egyptian fraction representations of~$r$ (that is, such that the sum of the reciprocals of the elements in the subset equals~$r$). If~$S$ is chosen to be of the form $\{1,2,\dots,n\}$, then the algorithm naturally finds Egyptian fraction representations with denominators bounded by $n$; indeed in this way we can determine the best possible bound (the minimal largest denominator) among Egyptian fraction representations of a given rational number~$r$. However, our algorithm works with an arbitrary finite set of positive integers, and so can be used to find Egyptian fraction representations with any desired restrictions on the denominators. (Indeed, the positive integers in the input set~$S$ don't even need to be distinct: the algorithm works perfectly well for multisets, although the resulting representations are not technically Egyptian fractions.) Our code is available online~\cite{Shi2017}.

\section{Motivation for the algorithm}

In this section, we set some notation to be used throughout the paper, and then expound the rationale behind our algorithm with a top-down approach. We give pseudocode for a preliminary version \textsc{ufrac-prelim} of our algorithm and work through some simple examples to orient the reader; later, we introduce the notion of ``reserved denominators'' to refine this preliminary algorithm. (A complete pseudocode description of the implemented algorithm \textsc{ufrac} will be given in Section~\ref{sec.algorithm} below.) At the end of this section we highlight a mathematical lemma that is crucial to the implementation.

\begin{notation}\label{notation}
Throughout this report, $r$ denotes a rational number (which in practice will be nonnegative), while~$D$ and~$D'$ denote finite multisets of positive integers (that is, ``sets'' in which elements are allowed to appear with any positive-integer multiplicity). We will always write rational numbers in reduced form, and thus we can unambiguously refer to the numerator and denominator of~$r$; for example, if $r$ is an integer then its denominator equals~$1$.

We define $R(D)$ to be the reciprocal sum of $D$ and $\delta(D,r)$ to be the difference between this reciprocal sum and the target rational~$r$, that is,
\begin{equation} \label{R and delta def}
R(D)=\sum_{d\in D} \frac1d \quad\text{and}\quad \delta (D,r) = R(D)-r
\end{equation}
(where the summand $\frac1d$ appears with the same multiplicity as $d$ does in $D$, if appropriate); by convention, $R(\varnothing)=0$.

We say that $D$ is a \emph{representation} of $r$ if $R(D)=r$ (or equivalently if $\delta(D,r)=0$). When $D$ consists of distinct integers, this is an Egyptian fraction representation, with the slight misuse of terminology in that we say (for example) that the set $D=\{2,3,6\}$ is a representation of $r=1$ when technically it is $R(D) = \frac12+\frac13+\frac16=1$ that is the Egyptian fraction representation.
\end{notation}

The following theorem describes inputs and the output of our algorithm.

\begin{theorem}
Let~$D$ be a multiset of positive integers, and let~$r$ be a rational number. The procedure~\textsc{UFRAC} defined in Section~\ref{sec.algorithm} generates the set of all submultisets of~$D$ that are representations of~$r$:
  \[\textsc{UFRAC}(D,r)=\bigl\{D'\subset D \colon R(D') = r\bigr\}.\]
\end{theorem}

For the remainder of this section, we use pseudocode and examples to motivate the algorithm that will be presented fully in Section~\ref{sec.algorithm}.

\subsection{A preliminary algorithm}

Our algorithm focuses on the difference $\delta (D,r) = R(D)-r$ (recalling that $D$ is a representation of~$r$ precisely when $\delta(D,r)=0$). If $\delta(D,r)<0$, it is impossible for a submultiset of $D$ to form a representation of~$r$, since the reciprocal sum of the submultiset is less than or equal to $R(D)$ and hence less than~$r$ (here we fundamentally use the positivity of the elements of~$D$). If $\delta(D,r)>0$, then it may be possible to find some representations of~$r$ after removing some carefully selected numbers from~$D$; this process is the main content of the algorithm.

Typically, $\delta(D,r)$ is not an integer, and thus its denominator is divisible by various powers of primes.
The idea is to remove some numbers in $D$ so that the greatest prime power $p^t$ dividing the denominator of $\delta(D,r)$ decreases, and then to repeat recursively. It is easy to check that the presence of $p^t$ in the denominator is determined exclusively by the multiples of $p^t$ in~$D$ (this is the phenomenon that explains, for instance, why one cannot find a representation of $r=\frac23$ using only denominators in $\{n\in\mathbb{N} \colon 3\nmid n\}$).
Hence it suffices to examine the multiples of $p^t$ when choosing numbers to remove. This strategy of focusing on one prime power at a time is central to the results of~\cite{Martin1999,Martin2000} and is also central to our algorithm.

Occasionally $\delta(D,r)$ is actually a positive integer. In this case, we identify the least number $\ell$ in~$D$ and create two branches (recursive calls to \textsc{ufrac-prelim}), one in which we retain~$\ell$ and insist that it be an element of any further representations, and one in which we remove~$\ell$. 

Using this strategy, we can give the following preliminary version of the algorithm.

\begin{algorithmic}[1]
  \Procedure{ufrac-prelim}{$D$, $r$}
  \If{$\delta(D,r)<0$}
  \State{\Return{$\varnothing$}}
  \ElsIf{$\delta(D,r)=0$}
  \State{\Return{$\{D\}$}}
  \ElsIf{$\delta(D,r)>0$}
  \If{$\delta(D,r)$ is an integer}
  \State{Let $\ell$ be the least element of $D$}
  \State{Let $D'$ be the submultiset of $D$ so that $D'\cup\{\ell\}=D$}
  \Statex{}\Comment{we write $D'=D\setminus\{\ell\}$, as we would for sets}
  \State{Let $R$ be the result of appending $\ell$ to each element of $\textsc{ufrac-prelim}(D',r-\frac1\ell)$}
  \State{\Return{$\textsc{ufrac-prelim}(D',r) \cup R$}}\label{efrac:line:branch}
  \Else{}
  \State{Let $p^{t}$ be the greatest prime power that divides the denominator of $\delta(D,r)$}
  \State{Determine sets of multiples of $p^t$ to remove from $D$ to obtain submultisets $D'$ such that $p^t$ does not divide the denominator of $\delta(D',r)$}\label{efrac:line:kill}
  \Statex{}\Comment{Section~\ref{the way to choose subsets} describes how to determine these sets}
  \State{\Return{the union of $\textsc{ufrac-prelim}(D',r)$ for each such~$D'$}}
  \EndIf{}
  \EndIf{}
  \EndProcedure{}
\end{algorithmic}
\smallskip

The following example traces through this preliminary algorithm in a concrete case.

\begin{example}\label{ex1}
Suppose $D=\{2,3,4,12\}$ and $r=\frac13$. A quick mental calculation reveals that $\{3\}$ and $\{4,12\}$ are the only two subsets of~$D$ that are representations of $\frac13$. We walk through how $\textsc{ufrac-prelim}(D,r)$ performs this calculation.
  
Initially, $R\bigl(\{2,3,4,12\}\bigr)=\frac12+\frac13+\frac14+\frac1{12}=\frac76$ and $\delta(D,r)=\frac76-\frac13=\frac56$, so that $3$~is the greatest prime power dividing the denominator of $\delta(D,r)$. The multiples of~$3$ in~$D$ are precisely $\{3,12\}$. We can check that removing $\{3\}$ from~$D$ results in $\delta\bigl(\{2,4,12\},\frac13\bigr)=\frac12$, while removing $\{12\}$ from~$D$ results in $\delta\bigl(\{2,3,4\},\frac13\bigr)=\frac34$; in both cases, $3$ no longer divides the denominator, so $\{2,4,12\}$ and $\{2,3,4\}$ are deemed acceptable choices for~$D'$ in line~\ref{efrac:line:kill} of the above procedure. On the other hand, removing $\{3,12\}$ from~$D$ results in $\delta\bigl(\{2,4\},\frac13\bigr)=\frac5{12}$, whose denominator is still divisible by~$3$; thus we do not allow $D'=\{3,12\}$ in line~\ref{efrac:line:kill}. Consequently, there are precisely two branches, one for each acceptable choice of~$D'$.

\begin{itemize}
\item In the first branch, where $\delta\bigl(\{2,4,12\},\frac13\bigr)=\frac12$, the largest prime power dividing the denominator is~$2$; in this case, all elements $\{2,4,12\}$ are multiples of~$2$. It turns out that there are precisely two subsets of $\{2,4,12\}$ that we can remove to yield a difference whose denominator is not divisible by~$2$. If we remove $\{2\}$, the result is  $\delta\bigl(\{4,12\},\frac13\bigr)=0$, revealing that $\{4,12\}$ is a representation of~$\frac13$. If we instead remove $\{2,4,12\}$, the result is $\delta\bigl(\varnothing,\frac13\bigr)=-\frac13$, which does not lead to a representation.
\item In the second branch, where $\delta\bigl(\{2,3,4\},\frac13\bigr)=\frac34$, the largest prime power dividing the denominator is~$4$. In this case, only $4$ is a multiple of~$4$, and removing $\{4\}$ yields $\delta\bigl(\{2,3\},\frac13\bigr)=\frac12$ (which is technically a new branch). Now the largest prime power dividing the denominator is~$2$, and removing the only remaining even denominator~$2$ yields $\delta\bigl(\{3\},\frac13\bigr)=0$, revealing that $\{3\}$ is a representation of~$\frac13$ as well.
\end{itemize}
\end{example}

\subsection{Reserving and removing denominators}\label{reserving and removing denominators}

Our actual implementation keeps track of certain additional data about previous branches, which we now describe.

Notice that in the first step of Example~\ref{ex1}, if we choose to remove $\{3\}$ and retain $\{12\}$, then $12$ appears in the final representation $\{4,12\}$; on the other hand, if we choose to remove $\{12\}$ and retain $\{3\}$, then $3$ remains in the final representation $\{3\}$. This is no coincidence: at each step, we decide whether to retain or remove each multiple of the greatest prime power in the denominator simultaneously. It proves beneficial to mark the retained denominators as ``reserved'' and keep track of these as denominators that are forced to appear in representations resulting from the given branch. When we do so, Example~\ref{ex1} is now modified as follows.

\begin{example}\label{ex2}
Again let $D=\{2,3,4,12\}$ and $r=\frac13$, so that $R\bigl(\{2,3,4,12\}\bigr)=\frac76$ and $\delta(D,\frac13)=\frac76-\frac13=\frac56$. As before, $3$~is the greatest prime power dividing the denominator, so we examine the denominators $\{3,12\}$ that are multiples of~$3$, and we may either remove $\{3\}$ and reserve $\{12\}$ or remove $\{12\}$ and reserve $\{3\}$.
\begin{itemize}
\item In the first branch, where we remove $\{3\}$, the reserved numbers are $\{12\}$ and the denominators to be examined are $\{2,4\}$. Since we have reserved $\{12\}$, the problem reduces to finding representations of $r_1=\frac13-R\bigl(\{12\}\bigr)=\frac14$ in $\{2,4\}$. Thus instead of computing $\delta\bigl(\{2,4,12\},\frac13\bigr)$, we need to compute $\delta\bigl(\{2,4\},\frac14\bigr)$. Notice that $\delta\bigl(\{2,4\},\frac14\bigr)=\frac12$, so we need to remove the prime power $2$ from the denominator. The only possible way is to retain $\{4\}$ and remove $\{2\}$; we find of course that $\delta\bigl(\{4\},\frac14\bigr)=0$. Hence $\{4\}$ is the only representation of $r_1=\frac14$ in $\{2,4\}$. We recover the representation of $r=\frac13$ in this branch by taking the union of reserved numbers and the representation of $r_1=\frac14$, that is, $\{12\}\cup\{4\}=\{4,12\}$.
\item In the second branch, where we remove $\{12\}$, the reserved numbers are $\{3\}$ and the denominators to be examined are $\{2,4\}$. Since we have reserved $\{3\}$, the problem reduces to finding representations of $r_1=\frac13-R\bigl(\{3\}\bigr)=0$ in $\{2,4\}$, and clearly $\varnothing$ is the only such representation. We then recover the representation of $\frac13$ in this branch by taking union of the reserved numbers and the representation of $r_2=0$, that is, $\{3\}\cup\varnothing=\{3\}$.
\end{itemize}
\end{example}

The extra work of reserving numbers may seem burdensome, but it has many crucial benefits:

\begin{itemize}
\item Without reserving numbers, a denominator will be a multiple of different prime powers, and thus would need to be examined multiple times. When we either reserve or remove a number, however, we only need to examine a denominator at most once in each branch, which accelerates the algorithm by removing redundancy.
\item Tracking reserved denominators provides a tree structure among all the branches and their descendents (subbranches). Subbranches of a branch will differ from each other in at least one reserved number, so they will not lead to the same representation. (The situation where~$D$ is a multiset needs attention in this regard, and we explicitly remove redundant branches caused by multiple occurrences of a given denominator.)
\item No matter how a subbranch is created, the number of unexamined denominators will be strictly less than the number of unexamined denominators in its parent branch. This reduction guarantees that the algorithm will terminate in a finite number of steps.
\item\label{reservation:3} The central idea in our algorithm is to treat one prime power divisor of the denominator of $\delta(D,r)$ at a time; however, this strategy needs to be assisted at any point when $\delta(D,r)$ is a positive integer. Reserving denominators provides a natural way to address this situation, by arbitrarily singling out the least number~$\ell$ in~$D$. We form one subbranch by removing (one copy of)~$\ell$, letting $D'=D\setminus\{\ell\}$ and continuing with $\delta(D',r)$. (As with any branch, this branch will be dropped if $\delta(D',r)<0$.) We form a second subbranch by instead reserving~$\ell$, so that $D'=D\setminus\{\ell\}$ still but we continue with $\delta(D',r')$ where $r'=r-\frac1\ell$; all representations of $r'$ using denominators in~$D'$ become representations of~$r$ with denominators in~$D$ when~$\ell$ is appended to them.
\end{itemize}

\begin{Rem}\label{remark reserving}
For the sake of later reference, it is helpful to explicitly record the relationships described in the last bullet point above. Let~$D$ be a multiset of positive integers (the unexamined denominators in a particular branch) and~$r$ a rational number (the target of that branch).
\begin{enumerate}
\item[(a)] If we reserve a submultiset~$E$ of denominators from~$D$, then the multiset of unexamined denominators in the resulting subbranch will be $D'=D\setminus E$, the target of the new branch will be $r'=r-R(E)$, and $\delta(D',r')=\delta(D,r)$.
\item[(b)] If we remove a submultiset~$E$ of denominators from~$D$, then the multiset of unexamined denominators in the resulting subbranch will be $D'=D\setminus E$, the target of the new branch will be $r'=r$, and $\delta(D',r')=\delta(D,r)-R(E)$.
\end{enumerate}
\end{Rem}

\subsection{A crucial number theory lemma}\label{the way to choose subsets}

We have seen that each branch targets the largest prime power dividing the denominator of $\delta(D,r)$ at that moment. So far, we have omitted any details about how to choose the submultisets of multiples of that prime power that cause it to no longer divide that denominator. (In Examples~\ref{ex1} and~\ref{ex2} we simply found the acceptable subsets by hand.) For the general situation, our procedure relies upon the following lemma, which is an easy exercise in elementary number theory.

\begin{lemma}\label{lemma}
Let $p$ be a prime and let $s\ge1$. Let $m/np^s$ be a positive rational number with $p\nmid n$, and let $c_1,\dots,c_k$ be integers not divisible by~$p$. For any subset $J$ of $\{1,\dots,k\}$, the denominator of $\displaystyle\frac m{np^s}-\sum_{j\in J} \frac1{c_j p^s}$ is not divisible by $p^s$ if and only if
  \begin{equation} \label{lemma-cond}
    mn^{-1}\equiv\sum_{j\in J}c_j^{-1}\pmod{p}.
  \end{equation}
\end{lemma}

Notice that the hypotheses of the lemma require that the denominators $c_j p^s$ not be divisible by $p^{s+1}$; however, this restriction might not hold for the set of multiples of $p^s$ in our denominator set~$D$. For instance, in the first branch of Example~\ref{ex1}, the greatest prime power dividing the denominator of $\delta\bigl(\{2,4,12\},\frac13\bigr)=\frac12$ is $2^1$; but we cannot apply the lemma to the denominators $\{2,4,12\}$ (if we tried to, we would see for instance that $12=6\cdot2^1$ and $6$ has no multiplicative inverse modulo~$2$). 

We address this issue in our algorithm by branching not on $p^t$, the greatest prime power dividing the denominator of $\delta(D,r)$, but rather on the greatest power $p^s$ dividing any of the elements of~$D$, with our goal being to remove the presence of~$p^s$ in subbranches. We can then choose integers $c_1,\dots,c_J$ not divisible by~$p$ such that the multiples of $p^s$ in~$D$ are precisely $\{c_1 p^s,\dots,c_J p^s\}$ and apply Lemma~\ref{lemma} to those integers. In every subbranch, each one of these $c_j p^s$ is either removed or reserved, and hence we no longer need to consider multiples of $p^s$ in subbranches. In particular, note that if $p^t<p^s$ and $p^s$ does not divide the denominator of~$r$, then removing all multiples of $p^s$ is a valid subbranch.

\section{Description of the full algorithm}\label{sec.algorithm}

We now fully describe (in pseudocode) our algorithm for finding all submultisets of a given input multiset whose reciprocal sum equals a target rational number. We need a data structure \emph{branch} to encase all the useful information at a given computational step of the procedure, a control structure \textsc{ufrac} managing the tree of computational branches, and a procedure \textsc{kill} to process branches and decide which new branches to make. We start with the data structure.

\newpage
\begin{Def}\label{branch def}
A \emph{branch} is a data structure that contains the following information:
\begin{itemize}
\item $D$, a multiset of positive integers;
\item \textit{rsvd}, a multiset of reserved numbers;
\item \textit{original-r}, the initial input $r$ to the procedure \textsc{ufrac};
\item $r$, the difference $\textit{original-r}-R(\textit{rsvd})$;
\item \textit{diff}, the rational number $\delta(D,r)$ (in reduced form)---recall the definitions of~$R$ and~$\delta$ from equation~\eqref{R and delta def};
\item \textit{gpp}, the greatest prime power dividing the denominator of \textit{diff}. (Note that~\textit{gpp}, which we denote by~$p^t$, is not necessarily the greatest power of $p$ dividing the elements of~$D$, which we denote by~$p^s$ and to which we will apply Lemma~\ref{lemma}.)
\end{itemize}
\end{Def}

\subsection{Reserving and removing denominators}

Remark~\ref{remark reserving} states the direct effects of removing and reserving denominators when generating new branches; however, the indirect effects need some special handling.
If we reserve a submultiset~$E$ of~$D$ of a branch $\textit{br}$, then $r$ will decrease by $R(E)$ as per Remark~\ref{remark reserving}~(a).
\begin{itemize}
\item If $r-R(E)<0$, then we drop this new branch. 
\item If $r-R(E)=0$, then the new branch is a representation. 
\item If $r-R(E)>0$, we then have to remove from $D$ all elements of $D\setminus E$ whose reciprocals are greater than $r-R(E)$ (since reserving such an element would cause the new $r$ to be negative).
\end{itemize}
Similarly, if we remove a submultiset~$E$ of~$D$ of $\textit{br}$, then $\textit{diff}$ will decrease by $R(E)$ as per Remark~\ref{remark reserving}~(b). 
\begin{itemize}
\item If $\textit{diff}-R(E)<0$, then we drop this new branch.
\item If $\textit{diff}-R(E)=0$, then the new branch is a representation. 
\item If $\textit{diff}-R(E)>0$, we will have to reserve all elements of $D\setminus E$ whose reciprocals are less than $\textit{diff}-R(E)$ (since removing such an element would cause the $\textit{diff}$ of a new branch to be negative).
\end{itemize}

The above observations show that the action of reserving elements in $D$ possibly needs to be followed by the action of removing elements in $D$ and vice versa. For this reason we have a procedure, which we call \textsc{reduce}, to correctly process all denominators that are greater than either $r$ or $\textit{diff}$.
\begin{algorithmic}[1]
  \Procedure{reduce}{\textit{br}}
  \While{there exists an element of $D$ that is greater than $r$ or \textit{diff}}
  \If{there exists an element of $D$ that is greater than $r$}
  \State{Create a new branch \textit{new-br} by removing all elements of $D$ greater than $r$}
  \Statex{}\Comment{as per Remark~\ref{remark reserving}~(b)}
  \Else{}\Comment{there exists an element of $D$ that is greater than \textit{diff}}
  \State{Create a new branch \textit{new-br} by reserving all elements in $D$ greater than \textit{diff}}
  \Statex{}\Comment{as per Remark~\ref{remark reserving}~(a)}
  \EndIf{}
  \State{Replace \textit{br} by \textit{new-br}}  
  \EndWhile{}
  \State{\Return{\textit{br}}}
  \EndProcedure{}
\end{algorithmic}

\noindent
From now on, we assume that all new branches are automatically processed by \textsc{reduce}.

\subsection{The control structure \textsc{ufrac}}

The input to \textsc{ufrac} is a multiset~$D$ of positive integers and a target nonnegative rational number~$r$ (both supplied by the user), and the output is all the representations of~$r$ using denominators in~$D$. We implement this procedure as a depth-first search to minimize the amount of memory required for the list of branches. The main mathematical work is performed by a subprocedure \textsc{kill}, described in the next section, that is called from \textsc{ufrac}.

\begin{algorithmic}[1]\label{ufrac}
  \Procedure{ufrac}{$D$, $r$}
  \State{Augment the data $D$ and $r$ into a complete branch \textit{br}}
  \If{\textit{diff} of \textit{br} is negative}
  \State{\Return{$\varnothing$}}
  \Else{}
  \State{Set \textit{branches} to be $\{\textit{br}\}$ and \textit{representations} to be $\varnothing$}
  \Statex{}\Comment{both \textit{branches} and \textit{representations} are ordered lists of branches}
  \While{\textit{branches} is not empty}
  \State{Set \textit{br} to be the first branch in \textit{branches}}
  \Statex{}\Comment{this \textit{br} is a local variable in the while-loop}
  \State{Set \textit{kill-result} to be the output of $\textsc{kill}(\textit{br})$}
  \Statex{}\Comment{this output is an ordered pair of ordered lists of branches}
  \State{Prepend the first ordered list of \textit{kill-result} to \textit{representations}}
  \State{Prepend the second ordered list of \textit{kill-result} to \textit{branches}}
  \Statex{}\Comment{prepending ensures that we are doing a depth-first search}
  \State{Remove \textit{br} from \textit{branches}}
  \EndWhile{}
  \State{Extract the multiset of denominators from each branch in \textit{representations}}
  \State{\Return{these multisets}}
  \EndIf{}
  \EndProcedure{}
\end{algorithmic}

\subsection{The Procedure \textsc{kill}}

As we saw in the above procedure, each branch is processed using a procedure \textsc{kill}, so named because its main purpose is to remove the largest prime power from the denominator of \textit{diff}. (Recall that $\textit{diff}=0$ indicates that we have found an exact representation of \textit{original-r}.) Because our treatment of the branch depends on both the sign and the integrality of \textit{diff}, we use \textsc{kill} to differentiate among the four relevant situations, introducing the procedures \textsc{kill-when-diff-is-integer} and \textsc{kill-when-diff-is-not-integer} which will be described shortly. The input to \textsc{kill} is a branch (as described in Definition~\ref{branch def}), while the output is a pair of collections of branches: the first collection contains representations of \textit{original-r} that happened to be found within this branch, and the second collection contains branches that have yet to be processed completely.

\newpage

\begin{algorithmic}[1]
  \Procedure{kill}{\textit{br}}
  \If{\textit{diff} of \textit{br} is $0$}
  \State{\Return{$(\{\textit{br}\},\varnothing)$}}
  \ElsIf{$\textit{diff}<0$ or $r<0$}\label{kill-cond}
  \State{\Return{$(\varnothing,\varnothing)$}}
  \ElsIf{$\textit{diff}$ is an integer}
  \State{\Return{$\textsc{kill-when-diff-is-integer}(\textit{br})$}}
  \Else{}\Comment{the case when $\textit{diff}$ is not an integer}
  \State{\Return{$\textsc{kill-when-diff-is-not-integer}(\textit{br})$}}
  \EndIf{}
  \EndProcedure{}
\end{algorithmic}


\subsubsection{The case when $\delta(D,r)$ is a positive integer}\label{kill when integer section}

It could happen that we have successfully removed all the primes dividing the denominator of \textit{diff}, but that $R(D)$ equals~$r$ plus some integer rather than~$r$ itself. The procedure \textsc{kill-when-diff-is-integer} handles this by choosing an element~$\ell$ from~$D$ (we choose the minimum element so that \textit{diff} is decreased most rapidly) and branching on whether~$\ell$ is or is not to be an element of the representations sought. This branching uses two new procedures, which are special cases of Remark~\ref{remark reserving}.
\begin{itemize}
\item The first procedure \textsc{br-reserve-l} generates a new branch, with~$\ell$ moved from~$D$ to \textit{rsvd} and~$r$ therefore replaced by $r-\frac1\ell$.
If $r<\frac1\ell$, we discard $\textsc{br-reserve-l}(\textit{br})$ as there will be no representation of $r$ that uses~$\ell$; if $r=\frac1\ell$, we promote $\textsc{br-reserve-l}(\textit{br})$ to a representation of~$r$ that has already been found.
\item The second procedure \textsc{br-remove-l} generates another new branch, with~$\ell$ simply removed from~$D$ leaving everything else unchanged. Since this branch will be processed by \textsc{kill} again, we don't need to check other conditions such as $R(D) \ge r + \frac1\ell$ (see line~\ref{kill-cond} of \textsc{kill} for example).
\end{itemize}
\smallskip

\begin{algorithmic}[1]
  \Procedure{kill-when-diff-is-integer}{\textit{br}}
  \State{Set $\ell$ to be the minimum number in D}
  \If{$r<\frac1\ell$}
  \State{\Return{$\bigl(\varnothing,\{\textsc{br-remove-l}(\textit{br})\}\bigr)$}}
  \ElsIf{$r=\frac1\ell$}
  \State{\Return{$\bigl(\{\textsc{br-reserve-l}(\textit{br})\},\{\textsc{br-remove-l}(\textit{br})\}\bigr)$}}
  \Else{}\Comment{this is the case when $r>\frac1\ell$}
  \State{\Return{$\bigl(\varnothing,\{\textsc{br-reserve-l}(\textit{br}),\textsc{br-remove-l}(\textit{br})\}\bigr)$}}
  \EndIf{}
  \EndProcedure{}
\end{algorithmic}

\subsubsection{The case when $\delta(D,r)$ is not an integer}

In this most common and most important step in our algorithm, we use Lemma~\ref{lemma} to remove elements from~$D$ to simplify the denominator of $\delta(D,r)$ in the following manner. We first compute the greatest prime power, $p^t$, dividing the denominator of $\delta(D,r)$; then, with that prime~$p$, we set~$s$ to be the largest integer such that $p^s$ divides at least one element of~$D$. We then collect all of the multiples of $p^s$ in~$D$, choosing integers $c_1,\dots,c_k$ so that
\[
M=\{d\in D \colon p^s\text{ divides }d\}=\{c_1p^s,c_2p^s,\dots c_k p^s\}.
\]
Note that~$p$ does not divide any of $c_1,\dots,c_k$ by the definition of~$s$.

We then use the procedure \textsc{kill-when-diff-is-not-integer} (whose pseudocode we will give after this discussion and two examples) to find all submultisets~$M'$ of~$M$ that satisfy
\begin{equation} \label{need s bigger than t}
\sum_{c_j p^s\in M'}c_j^{-1} \equiv \text{numerator of }\delta(D,r) \cdot p^{s-t} \cdot {\biggl(\frac{\text{denominator of }\delta(D,r)}{p^t}\biggr)}^{-1} \pmod{p}.
\end{equation}
(Each such submultiset $M'$ produces a new branch in our algorithm. We remark that if $p^s$ does not divide the denominator of~$r$, then $M'=M$ itself is always one valid choice.) By Lemma~\ref{lemma}, applied with $m = \bigl( \text{numerator of }\delta(D,r) \bigr) \cdot p^{s-t}$ and $n = \bigl( \text{denominator of }\delta(D,r) \bigr)/p^t$, removing $M'$ from $D$ results a multiset $D'$ such that $p^s$ does not divide the denominator of $\delta(D',r)$. Indeed, we both remove the elements of $M'$ from $D$ and also reserve the elements of $M\setminus M'$ in the branch so created.

The following example demonstrates this process.

\begin{example}\label{ex1714}
  Let $D=\{1,7,14,21,28\}$ and $r=1$, so that
  \[\delta(D,r)=\biggl( \frac11 + \frac17+\frac1{14}+\frac1{21}+\frac1{28} \biggr) - 1 = \frac{25}{3\times4\times7}.\]
Here $p^t=p^s=7^1$ and $M=\{7,14,21,28\}$. Next, we compute the right-hand side of the congruence~\eqref{need s bigger than t}, which is
  \[25\cdot1\cdot{(3\times4)}^{-1}\equiv5\pmod{7}.\]
  The multiplicative inverses of $1,2,3,4\pmod 7$ are $1,4,5,2$ respectively. Checking every possible subsum of $\{1,4,5,2\}$, we find that $1+4\equiv5\pmod{7}$, $5\equiv5\pmod{7}$ and $1+4+5+2\equiv5\pmod{7}$ (this last success was guaranteed since $7$ does not divide the denominator of~$r=1$). Therefore, we obtain new branches by removing $\{7,14\}$, $\{21\}$, and $\{7,14,21,28\}$, respectively, and reserving the remaining multiples of~$7$:
\smallskip
  \begin{center}
    \begin{tabular}{cccccc}
      \toprule
      \#& $D$ &\textit{rsvd}&$r$&\textit{original-r}&\textit{diff}\\
      \midrule
      1&$\{1\}$&$\{21,28\}$&$\frac{11}{12_{\mathstrut}}$&$1$&$\frac1{12}$\\
      2&$\{1\}$&$\{7,14,28\}$&$\frac3{4_{\mathstrut}}$&$1$&$\frac14$\\
      3&$\{1\}$&$\varnothing$&$1$&$1$&$0$\\
      \bottomrule
    \end{tabular}
  \end{center}
\smallskip
As desired, $7$ no longer divides the denominator of \textit{diff} in the resulting branches.

\end{example}

There is one possibility we have not yet considered: if $p^t > p^s$, then the right-hand side of the congruence~\eqref{need s bigger than t} is not well-defined due to the term $p^{s-t}$, and indeed in this situation there will never be any representations of~$r$ in~$D$.

\begin{example}
  Let $D=\{2,2,4,4\}$ and $r=\frac98$, so that $\delta(D,r)=\frac38$. Here, $p^t=8$ but $p^s=4$. Since there is no multiple of $8$ in~$D$, we cannot find a representation of~$\frac98$ in~$D$.
\end{example}

In our algorithm, we define a procedure \textsc{compute-pre-aim} that computes the prime powers~$p^t$ and~$p^s$ and, depending on their relative sizes, either returns the right-hand side of the congruence~\eqref{need s bigger than t} or returns a \textsc{false} flag to indicate that the branch can be abandoned.

This process requires a procedure \textsc{generate-subsets-aim} to find all of the submultisets that satisfy the congruence~\eqref{need s bigger than t}. We use an existing algorithm that solves the associated integer subset-sum problem~\cite{Cormen09}, which we will not explicitly describe here.

With this context in mind, we can now give the complete pseudocode for the procedure \textsc{kill-when-diff-is-not-integer}.


\begin{algorithmic}[1]
  \Procedure{kill-when-diff-is-not-integer}{\textit{br}}
  \State{Set \textit{pre-aim} to be $\textsc{compute-pre-aim}(\textit{br})$}
  \Statex{}\Comment{the quantities $p^t$ and $p^s$ are computed as by-products of \textsc{compute-pre-aim}}
  \If{\textit{pre-aim} is \textsc{false}}
  \Statex{}\Comment{the case when $p^s$ is impossible to remove from the denominator of \textit{diff}}
  \State{\Return{$(\varnothing,\varnothing)$}}
  \Else{}\Comment{the case when $p^s$ can be removed from the denominator of \textit{diff}}
  \Statex{}\Comment{here \textit{pre-aim} equals the right-hand side of the congruence~\eqref{need s bigger than t}}
  \State{Partition $D$ into $M$, the multiples of $p^s$, and $N$, the nonmultiples of $p^s$}
  \State{Set \textit{subsets} to be $\textsc{generate-subsets-aim}(M,\textit{pre-aim})$}
  \State{Set \textit{new-brs} to be the collection of complements of the elements of \textit{subsets}}
  \State{}\Comment{$\textit{new-brs}=\{M\setminus M' \colon M'\in\textit{subsets}\}$}
  \State{Augment each element \textit{newbr} of \textit{new-brs} into a complete branch}
  \Statex{}\Comment{as per Definition~\ref{branch def}; \textit{newbr} will be added to \textit{rsvd}, and $D$ will equal $N$}
  \State{Assign each element of \textit{new-brs} to \textit{new-representations} or \textit{new-branches}}
  \State{\Return{(\textit{new-representations}, \textit{new-branches})}}
  \EndIf{}
  \EndProcedure{}
\end{algorithmic}

\subsection{The early stopping version}

Sometimes we want to know merely whether there exists a representation of a number in a given multiset, in which case we do not need to compute all of the representations. We have therefore implemented an early stopping variant of the algorithm called \textsc{ufrac-early-stopping}. The input to \textsc{ufrac-early-stopping} is again a multiset~$D$ of positive integers and a target nonnegative rational number~$r$ (both supplied by the user), while the output is either one representation of~$r$ using denominators in~$D$, if any exist, or \textsc{false} if no such representations exist.

\begin{algorithmic}[1]
  \Procedure{ufrac-early-stopping}{$D$, $r$}
  \State{Make a branch \textit{br} from $D$ and $r$}
  \State{$\textsc{ufrac-early-stopping-recursion}(\{\textit{br}\})$}
  \EndProcedure{}
\end{algorithmic}
\smallskip

\noindent The procedure \textsc{ufrac-early-stopping-recursion} is structured as follows.

\newpage

\begin{algorithmic}[1]
  \Procedure{ufrac-early-stopping-recursion}{\textit{branches}}
  \If{\textit{branches} is empty}
  \State{\Return{$\varnothing$}}
  \Else{}
  \State{Take the first branch \textit{br} out from \textit{branches}}
  \State{Set \textit{kill-result} to be the output of $\textsc{kill}(\textit{br})$}
  \Statex{}\Comment{this output is an ordered pair of collections of branches}
  \State{Set \textit{new-representations} to be the first collection of \textit{kill-result}}
  \State{Set \textit{new-branches} to be the second collection of \textit{kill-result}}
  \If{$\textit{new-representations}$ is empty}
  \State{\Return{$\textsc{ufrac-early-stopping-recursion}(\textit{new-branches}\cup\textit{branches})$}}
  \Statex{}\Comment{prepending \textit{new-branches} ensures that we are doing a depth-first search}
  \Else{}
  \State{\Return{the multiset of denominators from the first element of \textit{new-representations}}}
  \EndIf{}
  \EndIf{}
  \EndProcedure{}
\end{algorithmic}


\section{Implementation and computational results}

As our algorithm involves many transformations of multisets, it is reasonable to use a functional programming language. We chose Scheme~\cite{Dybvig2009}, a dialect of Lisp, as our implementation language; our code is available online~\cite{Shi2017}.
Scheme has good native support for rationals with arbitrarily large numerators and denominators, limited only by the computer's memory; this functionality is not natively supported by many mainstream programming languages like C, C++, Java, or Python. Among all the implementations of Scheme, we chose Chez Scheme because it is the most efficient one in most cases~\cite{Schemebenchmarks}.

The main procedure, \textsc{ufrac}, takes a multiset~$D$ of positive integers and a target rational number~$r$ as input and returns all representations of~$r$ using denominators in~$D$. For example, the command
{\tt(ufrac (list 2 2 3 3 4 5 6 6 7 8 12) 3/2)}
yields a list of the five unit fraction representations of~$\frac32$ using the indicated denominators, namely
\smallskip
\centerline{\tt((2 3 4 6 6 12) (2 3 3 6 6) (2 3 3 4 12) (2 2 3 6) (2 2 4 6 12)).}
\smallskip
\noindent
In addition to entering~$D$ explicitly as above, one can also use any of Scheme's built-in functions or self-defined functions to construct~$D$: for example, the command
{\tt(ufrac (range 1 10) 3/2)}
yields a list of the two Egyptian fraction representations of~$\frac32$ using only denominators up to~$10$, namely
{\tt((1 2) (1 3 6))} (since $\frac11+\frac12=\frac32=\frac11+\frac13+\frac16$).

For another example, we mention a theorem of Graham~\cite{Graham1964} that a rational number~$r$ can be written as an Egyptian fraction using only denominators that are perfect squares if and only if $r \in \bigl[ 0, \frac{\pi^2}6-1 \bigr) \cup \bigl[ 1,\frac{\pi^2}6 \bigr)$. The Scheme command {\tt(map square (range 1 n))} produces a list of the first $n$ perfect squares. Using our code, the command {\tt(ufrac (map square (range 1 34)) 1/2)} yields {\tt\#f}, while the command {\tt(ufrac (map square (range 1 35)) 1/2)} yields {\tt((4 9 16 25 49 144 225 400 784 1225))}; this shows that among Egyptian fraction representations of $\frac12$ using only square denominators, the unique one with the smallest maximal denominator is $\tfrac12=\tfrac1{2^2}+\tfrac1{3^2}+\tfrac1{4^2}+\tfrac1{5^2}+\tfrac1{7^2}+\tfrac1{12^2}+\tfrac1{15^2}+\tfrac1{20^2}+\tfrac1{28^2}+\tfrac1{35^2}$.

Further examples and syntax can be found at~\cite{Shi2017}. We now describe two collections of computations we carried out that relate to extant open problems on Egyptian fractions.

\subsection{Smallest largest denominators}

We used our algorithm to explore the densest representations of integers, by which we mean the Egyptian fractions whose maximal denominators are as small as possible. Recalling that $R(D)=\sum_{d\in D} \frac{1}{d}$, we define
\[
G(r)=\min\bigl\{k\in\mathbb{N} \colon \text{there exists } D\subset\{1,2,\dots,k\} \text{ such that } R(D)=r\bigr\}
\]
to be the largest denominator of the densest Egyptian fraction representation of~$r$.
For example, $G(2)=6$, since $2=1+\frac12+\frac13+\frac16$ while there is no representation of~$2$ in $\{1,2,3,4,5\}$. We call the representation $\{1,2,3,6\}$ a \emph{witness} for the value $G(2)=6$.

We can compute $G(r)$ by simply executing the command {\tt(ufrac (range 1 n) r)} for increasing values of~$n$ until a representation is found. For example, the command {\tt(ufrac (range 1 23) 3)} returns {\tt\#f} while the command {\tt(ufrac (range 1 24) 3)} returns
\smallskip

\centerline{{\tt((1 2 3 4 5 6 8 9 10 15 18 20 24))},}
\smallskip
\noindent
which shows that $G(3)=24$ and that the witness is unique.

In this way, we found that $G(4)=65$, and indeed that
\begin{align*}
  4=R \smash{\bigl(} \{&1, 2, 3, 4, 5, 6, 7, 8, 9, 10, 11, 12, 13, 14, 15, 16, 18, 20, 22, 24,\\
            &26, 27, 28, 30, 33, 35, 36, 40, 42, 45, 48, 52, 54, 56, 60, 63, 65\} \smash{\bigr)}
\end{align*}
and that this witness is unique. We also found that $G(5)=184$; there are~$16$ witnesses, one of which gives

{\smaller\smaller\smaller
\begin{align*}
  5=R \smash{\bigl(} \{
  & 1,2,3,4,5,6,7,8,9,10,11,12,13,14,15,16,17,18,19,20,21,22,23,24,25,26,27,28,29,30,32,33, \\
  & 34,35,36,38,39,40,42,44,45,46,48,50,51,52,54,55,56,60,62,63,65,66,68,70,72,75,76,77,78, \\
  & 80,81,84,85,88,90,91,92,93,95,99,102,104,105,108,110,112,114,115,116,117,120,126,130,132, \\
  & 133,136,140,143,144,145,150,152,153,154,155,156,160,161,162,168,170,171,175,180,184
    \} \smash{\bigr)} .
\end{align*}
}

\noindent
We observed that~$136$ is an element of all~$16$ witnesses, a detail we will return to in a moment.

Finally, we found that $G(6)=469$; there are~$224$ witnesses, one of which gives

{\smaller\smaller\smaller
  \begin{align*}
  6=R \smash{\bigl(} \{
  & 1,2,3,4,5,6,7,8,9,10,11,12,13,14,15,16,17,18,19,20,21,22,23,24,25,26,27,28,29,30,31,32,33,34, \\
  & 35,36,37,38,39,40,41,42,43,44,45,46,47,48,49,50,51,52,53,54,55,56,57,58,60,61,62,63,64,65, \\
  & 66,67,68,69,70,72,74,75,76,77,78,80,81,82,84,85,86,87,88,90,91,92,93,94,95,96,98,99,100,102, \\
  & 104,105,106,108,110,111,112,114,115,116,117,118,119,120,121,122,123,124,126,128,129,130, \\
  & 132,133,134,135,138,140,141,143,144,145,147,148,150,152,153,154,155,156,159,160,161,162, \\
  & 164,165,168,170,171,174,175,176,180,182,183,184,185,186,187,188,189,190,192,195,196,198, \\
  & 200,201,203,204,205,207,208,209,210,212,215,216,217,220,221,224,225,228,230,231,232,234, \\
  & 238,240,242,245,246,247,248,250,252,253,255,258,259,260,261,264,266,268,270,272,273,275, \\
  & 276,280,282,285,286,287,288,290,294,295,296,297,299,300,301,304,305,306,308,310,312,315, \\
  & 319,320,322,324,325,328,329,330,333,336,340,341,342,344,345,348,350,351,352,357,360,363, \\
  & 364,368,372,374,375,376,377,378,380,384,385,387,390,396,399,400,402,405,406,408,413,414, \\
  & 416,418,420,424,425,429,430,432,434,435,440,442,444,448,451,455,456,460,462,465,468,469
    \} \smash{\bigr)} .
\end{align*}
}

\noindent
Using a recent laptop, the time it took to:
\begin{itemize}
\item verify that there are no representations of $r=6$ in $\{1,\dots,468\}$ was about $2.7$~s;
\item find one representation for $r=6$ in $\{1,\dots,469\}$ (using the early-stopping version of \textsc{ufrac}) was about $69$~s;
\item find all representations for $r=6$ in $\{1,\dots,469\}$ was about $560$~s.
\end{itemize}

These computations comprise independent verification of the first six entries of sequence A101877 in the Online Encyclopedia of Integer Sequences~\cite{OEIS}. It turns that that the next two values $G(7)=1243$ and $G(8) = 3231$ have been found by Hugo van der Sanden, with code that is available at~\cite{HvdS}. We were only able to verify that $G(7)>1210$. Searching for representations of~$7$ in $\{1,\dots,1230\}$, for example, ran for 56 hours and examined 21 million branches before we terminated the session manually.

We remark that~\cite{OEIS} reports the following question:

{\quote``Paul Hanna asks if it is always true that a solution set $S$ for $n+1$ must necessarily contain a solution set for $n$ as a subset. This is true for small $n$, apparently, but seems to me unlikely to hold in general.---N.~J.~A.~Sloane, Dec 31 2005''

}

\smallskip\noindent
We can answer this question in the negative: the witness given above for~$6$ does not contain~$136$, which means that it does not contain any witness for~$5$ since all of those witnesses do contain~$136$.

\subsection{Second-smallest denominators}

It is easy to see (using Lemma~\ref{lemma}, for example) that a prime cannot be the largest denominator in an Egyptian fraction representation of~$1$. Erd\"os asked about the density of the set of integers that cannot be the largest denominator in an Egyptian fraction representation of~$1$, and went on to ask the analogous question for the set of integers that cannot be the second-largest denominator in such a representation (and third-largest and so on). Perhaps unexpectedly, the first author~\cite[Theorem~2]{Martin2000} proved that every sufficiently large integer \textit{can} be the second-largest denominator in such a representation. We have carried out computations that support a speculation in that paper, which we promote to an explicit conjecture:

\begin{Conj}\label{2nd}
  Every integer $d\ge5$ can be the second-largest denominator in an Egyptian fraction representation of~$1$.
\end{Conj}

\noindent (It would follow easily from this conjecture and the splitting identity $\frac1n=\frac1{n+1}+\frac1{n(n+1)}$ that every integer $d\ge2$ can be the third-largest, fourth-largest, etc.~denominator in such a representation.) We have verified this conjecture for all $5\le d\le6{,}000$. The witnesses can be found at~\cite[{\tt conjecture\_witnesses.md}]{Shi2017}, and the first several are listed below; the final witness is the representation $\{5,6,7,8,9,14,15,19,38,56,95,114,6{,}000,2{,}394{,}000\}$ of~$1$, whose second-largest denominator is $d=6{,}000$.

\begin{table}[ht]
  \centering
  \begin{tabular}{ccc}
    \toprule
    $d$ & $c$ & representation \\
    \midrule
    5 & 4 & \{2, 4, 5, 20\} \\
    6 & 2 & \{2, 4, 6, 12\} \\
    7 & 6 & \{2, 3, 7, 42\} \\
    8 & 3 & \{2, 3, 8, 24\} \\
    9 & 2 & \{2, 3, 9, 18\} \\
    10 & 3 & \{2, 5, 6, 10, 30\} \\
    11 & 120 & \{3, 4, 5, 8, 11, 1{,}320\} \\
    12 & 2 & \{2, 4, 8, 12, 24\} \\
    13 & 12 & \{2, 4, 6, 13, 156\} \\
    14 & 6 & \{2, 4, 6, 14, 84\} \\
    15 & 4 & \{2, 4, 6, 15, 60\} \\
    16 & 3 & \{2, 4, 6, 16, 48\} \\
    17 & 16 & \{2, 4, 8, 16, 17, 272\} \\
    \bottomrule
  \end{tabular}
  \qquad
  \begin{tabular}{ccc}
    \toprule
    $d$ & $c$ & representation \\
    \midrule
    18 & 2 & \{2, 4, 6, 18, 36\} \\
    19 & 56 & \{2, 4, 8, 14, 19, 1{,}064\} \\
    20 & 3 & \{3, 4, 5, 12, 15, 20, 60\} \\
    21 & 20 & \{2, 4, 5, 21, 420\} \\
    22 & 10 & \{2, 4, 5, 22, 220\} \\
    23 & 45 & \{2, 3, 15, 18, 23, 1{,}035\} \\
    24 & 5 & \{2, 4, 5, 24, 120\} \\
    25 & 4 & \{2, 4, 5, 25, 100\} \\
    26 & 12 & \{2, 3, 8, 26, 312\} \\
    27 & 35 & \{2, 5, 7, 14, 21, 27, 945\} \\
    28 & 3 & \{2, 3, 14, 21, 28, 84\} \\
    29 & 28 & \{2, 4, 7, 14, 29, 812\} \\
    30 & 2 & \{2, 4, 5, 30, 60\} \\
    \bottomrule
  \end{tabular}
\end{table}

To verify Conjecture~\ref{2nd} for a particular integer $d\ge5$, we search for representations of~$1$ whose largest two denominators are $d$ and $c\times d$ for some well-chosen integer $c>1$, as indicated in the table. (When~$d$ is prime, it is easy to see that the largest denominator must necessarily be a multiple of~$d$; for simplicity we restrict to this form for all~$d$.) Among the integers~$c$ below a particular bound, we choose the one that minimizes the largest prime power factor of the denominator of $r=1-\frac{1}{d}-\frac{1}{cd}$. We then used our algorithm to search for representations of this~$r$ using only denominators up to~$100$ to limit runtimes (each such search took less than a tenth of a second on a laptop). In some cases, it was necessary to iteratively increase the bound~$100$ until a representation was found.

While this strategy occasionally needed to be adjusted by hand for small~$d$, it soon became clear that representations of the type we sought were plentiful, even under these artificial restrictions. The largest source of variation was that increasing the range from which we chose~$c$ made it more likely to find suitable representations of~$1$, while increasing that range made choosing~$c$ take longer: searching among $2\le c\le1{,}000$ takes only a few milliseconds, while searching among $2\le c\le500{,}000$ (which was needed as~$d$ increased) takes about $30$~seconds for each~$d$.

\section{Future work}

The algorithm has been brought to a stable and satisfying state. In the future, however, it would be attractive to build a multi-thread version of the algorithm to take advantage of parallel processing. We also remark that it seems difficult to determine the theoretical complexity of our algorithm; certainly we do not individually examine every possible submultiset of the input multiset, but some sort of exhaustive search is being performed in the end since the algorithm does produce every representation. Finally, we hope that others will find uses for the algorithm in their own Egyptian fraction investigations.

\providecommand{\bysame}{\leavevmode\hbox to3em{\hrulefill}\thinspace}
\providecommand{\MR}{\relax\ifhmode\unskip\space\fi MR }
\providecommand{\MRhref}[2]{%
  \href{http://www.ams.org/mathscinet-getitem?mr=#1}{#2}
}
\providecommand{\href}[2]{#2}


\begin{thebibliography}{99}

\bibitem{Bleicher2000}
A. Beck, M.~N. Bleicher, and D.~W. Crowe, \emph{Excursions into
  mathematics}, millennium ed., A K Peters, Ltd., Natick, MA, 2000, with a
  foreword by Martin Gardner. \MR{1744676}

\bibitem{Beeckmans1994}
L. Beeckmans, \emph{The splitting algorithm for Egyptian fractions},
J. Number Theory \textbf{43} (1993), no.~2, 173--185. 

\bibitem{Bleicher1972}
M.~N. Bleicher, \emph{A new algorithm for the expansion of {E}gyptian
  fractions}, J. Number Theory \textbf{4} (1972), 342--382. \MR{0323696}

\bibitem{Cormen09}
T. H. Cormen, C. E. Leiserson, R. L. Rivest, and C.~Stein, \emph{Introduction to
  algorithms}, 4 ed., The MIT Press, 2009.

\bibitem{Dybvig2009}
R. K. Dybvig, \emph{The {S}cheme programming language}, The MIT Press, The MIT
  Press, 2009.

\bibitem{Eppstein1995}
D. Eppstein, \emph{Ten Algorithms for Egyptian Fractions},
Mathematica in Education and Research \textbf{4} (1995), no.~2, 5--15. Available at \url{https://library.wolfram.com/infocenter/Articles/2926}; updated content at \url{https://www.ics.uci.edu/~eppstein/numth/egypt/binary.html}
  
\bibitem{GN2013}
E. Gyimesi and G. Nyul,
\emph{A note on Golomb's method and the continued fraction method for Egyptian fractions},
Ann. Math. Inform. \textbf{42} (2013), 129--134.

\bibitem{Graham1964}
R. L. Graham, \emph{On finite sums of unit fractions}, Proc. London Math. Soc. (3) \textbf{14} (1964),
193--207. MR{28:3968}

\bibitem{Martin1999}
G. Martin, \emph{Dense {E}gyptian fractions}, Trans. Amer. Math. Soc.
  \textbf{351} (1999), no.~9, 3641--3657. \MR{1608486}

\bibitem{Martin2000}
G. Martin, \emph{Denser {E}gyptian fractions}, Acta Arith. \textbf{95} (2000),
  no.~3, 231--260. \MR{1793163}

\bibitem{OEIS}
OEIS Foundation Inc., \emph{The On-Line Encyclopedia of Integer Sequences}. \url{http://oeis.org/A101877}

\bibitem{Rhindpapyrus}
\emph{The {R}hind {M}athematical {P}apyrus}, Amer. Math. Monthly \textbf{34}
  (1927), no.~9, 445--446. \MR{1521286}

\bibitem{Schemebenchmarks}
\emph{Scheme benchmarks}. \url{https://ecraven.github.io/r7rs-benchmarks/}

\bibitem{Shi2017}
Y. Shi, \emph{Dense Egyptian fractions}, GitHub repository.
  \url{https://github.com/ChlorophyII/Dense-Egyptian-Fractions}

\bibitem{Sylvester1880}
J.~J. Sylvester, \emph{On a point in the theory of vulgar fractions},
  Amer. J. Math. \textbf{3} (1880), no.~4, 332--335. \MR{1505274}
  
\bibitem{Takenouchi1921}
T. Takenouchi, \emph{On an indeterminate equation},
Proc. Physico-Mathematical Soc.\ of Japan (3rd ser.) \textbf{3}, 1921, 78--92.

\bibitem{HvdS}
H. van der Sanden, \emph{A101877}, GitHub repository. \url{https://github.com/hvds/seq/tree/master/A101877}
  
\bibitem{Wagon2010}
S. Wagon, \emph{Mathematica${}^{\text{\textregistered}}$ in Action: Problem solving through visualization and computation}, 3rd ed., Springer-Verlag New York, 2010.

\end{thebibliography}
\end{document}